\documentclass[12pt]{article}
\usepackage{amsfonts}
\usepackage{amsmath}
\usepackage{amsxtra}
\usepackage{amssymb,latexsym}
\usepackage[mathcal]{eucal}
\usepackage{amscd}

\newtheorem{theo}{{\bfseries Theorem}}[section]
\newtheorem{prop}[theo]{{\bfseries Proposition}}
\newtheorem{lem}[theo]{{\bfseries Lemma}}
\newtheorem{cor}[theo]{{\bfseries Corollary}}
\newtheorem{df}[theo]{{\bfseries Definition}}

\def \N {\mathbb N}
\def \C {\mathcal C}
\def \Z {\mathbb Z}
\def \R {\mathbb R}
\def \A {\mathcal A}
\def \B {\mathcal B}

\def \ep {\epsilon}

\def \d {\delta}
\def \l {\lambda}

\def \o {\omega}

\def \proof { {\bfseries Proof:} }

\usepackage{amssymb,latexsym}
\usepackage[mathcal]{eucal}
\usepackage{amscd}
\usepackage{amsmath}
\usepackage{graphicx}
\usepackage{amsfonts}
\usepackage{amscd}

\numberwithin{equation}{section}

\begin{document}

\title{\bfseries  Dynamics of Induced Systems}
\vspace{1cm}
\author{Ethan Akin, Joseph Auslander, Anima Nagar\\ \\
%\address
Mathematics Department,\\
 The City College, 137 Street and Convent Avenue,\\
 New York City, NY 10031, USA\\
%\email
ethanakin@earthlink.net\\ \\
%\address
Mathematics Department,\\
University of Maryland,\\
College Park, MD 20742, USA\\
%\email
jna@math.umd.edu \\ \\
%\address
Department of Mathematics,\\
Indian Institute of Technology Delhi,\\
Hauz Khas, New Delhi 110016, INDIA\\
%\email
anima@maths.iitd.ac.in}

    \vspace{.2cm}
\date{December, 2014}

\vspace{.2cm} \maketitle

\begin{abstract}

In this paper, we study the  dynamical properties of actions on the space of compact subsets of the phase space. More precisely, if $X$ is a metric space, let $2^X$ denote the space of non-empty compact subsets of $X$ provided with the Hausdorff topology. If $f$ is a continuous self-map on $X$, there is a naturally induced continuous self-map $f_*$ on $2^X$. Our main theme is the interrelation between the dynamics of $f$ and $f_*$. For such a study, it is useful to consider the space $\C(K,X)$ of continuous maps from a Cantor set $K$ to $X$ provided with the topology of uniform convergence, and $f_*$ induced on $\C(K,X)$  by composition of maps. We mainly study the properties of transitive points of the induced system $(2^X,f_*)$ both topologically and dynamically, and give some examples. We also look into some more properties of the system $(2^X,f_*)$.

\end{abstract}

\section{Introduction}

We first review some useful results and set up notation. All our basic definitions and notations are as in Akin \cite{ae}, though we make full attempt to explain each of these as we proceed.  All our spaces are assumed to be nonempty.

Recall that any two perfect, compact, zero-dimensional metric spaces are homeomorphic and we will refer to any such as a Cantor
set. A Polish space is a  separable space which admits a complete metric. As a compact metric space is complete and separable, it is
Polish.  A Polish space with no isolated points is called perfect.
Any perfect Polish space contains a Cantor set.  A $G_{\delta}$ subset of a Polish space is Polish. It follows that
any nonempty open subset of a perfect Polish space contains a Cantor set. For an exposition of this material, see, e.g. Akin \cite{aecm}.

Every compact metric space is a  continuous image of a Cantor set.  We include the brief proof.

\begin{prop}\label{prop1.1} If $X$ is a nonempty, compact, metrizable space and $K$ is a Cantor set then there is a continuous map from $K$ onto
$X$. \end{prop}

\proof Let $\mathcal B = \{ U_1, U_2,... \}$ be a countable base for $X$ and let $Q $ be the closure in $X \times  \{ 0,1 \}^{\N}$
of $\{ (x,a) : a_i = 1 \Leftrightarrow x \in U_i \}$. Clearly, $Q$ is contained in the closed set
$\{ (x,a) : a_i = 1 $ and $x \in \bar U_i$ or $a_i = 0 $ and $x \in X \setminus U_i \}$.
It follows that if $(x_1,a_1), (x_2,a_2) \in Q$
and $U_i \in \B$  with $x_1 \in U_i, x_2 \in X \setminus U_i$ then $(a_1)_i = 1, (a_2)_i = 0$. Hence, the projection map from
$Q$ to $\{ 0,1 \}^{\N}$ is injective and so $Q$ is a zero-dimensional compact metric space.
Clearly, $\pi_1 : Q \to X$ is surjective. $Q \times K$ is perfect and so there is a homeomorphism $h :K \to Q \times K$.
$\pi_1 \circ h : K \to X$ is surjective.

$\Box$ \vspace{.5cm}

For a metric space $X$, we follow Illanes and Nadler \cite{in} letting
$2^X$ denote the space of nonempty compact subsets of $X$.  Notice that we are excluding
$\emptyset \subset X$ which is sometimes regarded as an isolated point of $2^X$. If $K$ is a compact metrizable space
we let $\C(K,X)$ denote the space of continuous maps from $K$ to $X$.  Each of these has a natural metric induced from $d$ on $X$.
(We will use $d$ for all the metrics which occur, allowing context to make the referent clear.)

On $2^X$ we use the Hausdorff metric:  For $A, B \in 2^X$
\begin{equation}\label{1.1}
d(A,B) \quad = \quad max \{d(a,B) : a \in A  \} \cup \{ d(b,A) : b \in B  \},
\end{equation}
where $d(a,B) = min \{ d(a,b) : b \in B \}$. Thus, $d(A,B) < \ep $ if and only if each set is in the open $\ep $ neighborhood of the other,
or, equivalently, each point of $A$ is within $\ep$ of a point in $B$ and vice-versa.

When $X$ is compact, we occasionally use an equivalent topology on $2^X$. Define for any collection $\{ U_i : 1 \leq i \leq n \}$ of opene ( = open and nonempty) sets,
\begin{equation}
<U_1, U_2, \ldots U_n> =  \{E \in 2^X :E \subseteq \bigcup \limits_{i=1}^n
U_{i}, \ E \bigcap U_{i} \neq \phi,  \textrm{ } 1 \leq i \leq n \}
 \end{equation}
The topology on $2^X$, generated by such collection as basis, is
known as the {\emph{Vietoris topology}}.

As topological objects, $2^X$ have a very rich structure. It is known that if $X$ contains a \emph{Peano continuum} then $2^X$ contains the Hilbert Cube. We refer  the reader to Illanes and Nadler \cite{in}  or Schori and West \cite{sw} for more details.

\vskip .1cm

On $\C(K,X)$ we use the sup metric: For $u,v \in \C(K,X)$
\begin{equation}\label{1.2}
d(u,v) \quad = \quad \max \{ d(u(x),v(x)) : x \in K \},
\end{equation}
with the topology of uniform convergence.

If $\{ A_n \}$ is a sequence of closed sets in a topological space then
\begin{equation}\label{1.3}
\begin{split}
\overline{ \bigcup_n \{ A_n \}} \quad = \quad \bigcup_n \{ A_n \} \ \cup \ Lim sup_n \{ A_n \}, \hspace{1cm}\\
\mbox{where} \qquad \ Lim sup_n \{ A_n \} \quad = \quad \ \bigcap_k \overline{\bigcup_{n \geq k} \{ A_n \}}.
\end{split}
\end{equation}

The results stated below are standard. One can refer to \cite{ae, in, mi} for more details. For the sake of completeness we briefly discuss the proofs of the results which we shall use and in the process establish some of our notations.

\begin{lem}\label{lem1.2} For a metric space $X$, and $ \{ A_n \}$  a sequence in $2^X$,

(a) If $\{ A_n \}$ converges to $A$ in $2^X$ then $\overline{ \bigcup_n \{ A_n \}}$ is compact.

(b) If $\overline{ \bigcup_n \{ A_n \}}$ is compact and $\{ A_n \}$ is Cauchy then $A_n$ converges to  $ Lim sup \{ A_n \}$.

(c) If $X$ is complete then $2^X$ is complete.

(d) If $X$ is compact then $2^X$ is compact.\end{lem}

\proof We consider the sequence $ \{ A_n \}$  in $2^X$ for a metric space $X$.

(a) Let $\{x_n \}$ be a sequence in $\overline{ \bigcup_n \{ A_n \} }$. If there is an $N \in \N$ such that infinitely many $x_i$ are in the compact set $\cup _{n \leq N} A_n$ then clearly $\{x_n\}$ has a convergent subsequence. Otherwise, since $A_n \to A$, it follows from the definition of the Hausdorff metric that there are $y_i$ in $A$ with $d(x_i,y_i) \to 0$. Since $A$ is compact, $\{y_i\}$ has a convergent subsequence, and the corresponding subsequence of $\{x_i\}$ converges to the same point.

(b) Let $\epsilon > 0$. Since $\{A_n\}$ is a Cauchy sequence there is an $N>0$ such that all the $A_n$ for $n\geq N$ are within $\epsilon/3$ neighbourhood of each other, and so within $2 \epsilon/3$ of $\overline{ \bigcup_n \{ A_n \} }$.
Now $\overline{\cup _{k \geq N} A_k} \to Lim sup A_n$ as $ k \to \infty$, so $A_n$ is within $\epsilon$  neighbourhood of $Lim sup A_n$  when $n \geq N$. It follows that $\lim A_n = Lim sup A_n$. In particular if $A_n \to A$ then $A= Lim sup A_n$.

(c) Let $\{A_n\}$ be a Cauchy sequence in $2^X$. We show that $\overline{\cup A_n}$ is compact. Since a compact metric space is complete, this will imply that $\{A_n\}$ converges. Now since $\overline{\cup A_n}$ is closed, it is sufficient to show that it is totally bounded. Let $\epsilon > 0$. Again there is an $N \in \N$ such that all $A_n$ for $n \geq N$ are within $\ep/3$ neighbourhood of $A_N$. Therefore a finite subset of $A_N$ which is $\epsilon/3$ dense in $A_N$ is $\epsilon$ dense in $\overline{\cup _{n \geq N} A_k}$. Since $\cup_{n<N} A_n$ is compact, there is a finite $\epsilon$ dense set in $\overline{\cup  A_n}$.

(d) If $X$ is compact, and $\ep > 0$ then there exists $D$  a finite subset of $X$ which is
$\ep$-dense, i.e. for every $x \in X$ there exists $d \in D$ such that $d(x,d) < \ep$. If $A \in 2^X$ then the set of points
of $D$ which are less than $\ep$ from a point of $A$ form an element of $2^X$ which has distance less than $\ep$ from $A$.
Thus, the dense set of finite subsets of $2^X$ are also $\ep$-dense. Thus $2^X$ is  totally bounded as well as complete, and so  is
compact.

$\Box$ \vspace{.5cm}

 We will let $(n)$ denote $\{ 1,...,n \}$ and write
$X^{(n)}$ for the $n$-fold product of copies of $X$. Hence, $f^{(n)}$ denotes the $n$-fold product of copies of $f$ as
opposed to $f^n$ which is the $n$-fold iterate of a map $f$ on $X$.

\begin{prop}\label{prop1.3} Let $X$ be a metric space and $K$ be a compact metrizable space.
\begin{itemize}
\item[(i)]  Define $i_X : X \to 2^X$ by $x \mapsto \{ x \}$ and
$i_X : X \to \C(K,X)$ by $x \mapsto c_x$ where $c_x$ is the function
which takes the constant value $x$. Both maps
% eja del$i_X$
are isometric inclusions.

\item[(ii)] Define $i_{X,n} : X^{(n)} \to 2^X$ by $(x_1,...,x_n) \mapsto \{ x_1,..., x_n \}$ and $\vee : 2^{2^X} \to 2^X$ by
$F \mapsto \bigcup F$. Each of these has Lipschitz constant $1$, where on the product space $X^{(n)}$ we use the max of the
distances between corresponding coordinates. Hence, the map $(2^X)^{(n)} \to 2^X$
given by $(A_1,...,A_n) \mapsto A_1 \cup ... \cup A_n$
has Lipschitz constant $1$. The image of $i_{X,n}$ is closed and $\vee$ is surjective.

\item[(iii)] For $D \subset X$ let $FIN(D) \subset 2^X$ be the collection of all finite subsets of $D$. For $\A$ a clopen partition
of $K$ we let $\C(\A,D) \subset \C(K,X)$ denote the functions with values in $D$ which are constant on each element of $\A$.
If $D$ is dense in $X$ then $FIN(D)$ is dense in $2^X$. If $K$ is a Cantor set and  $\{ \A_n \}$ is a sequence of clopen partitions
whose mesh tends to zero then $\bigcup_n \{ \C(\A_n,D) \}$ is dense in $\C(K,X)$. In particular, if $X$ is separable then
so are $\C(K,X)$ and $2^X$. Furthermore, if $X$ has no isolated points then neither does $2^X$ and if $K$ is a Cantor set
then $\C(K,X)$ has no isolated points either.

\item[(iv)] The map $Im : \C(K,X) \to 2^X$ defined by $Im(u) = u(K)$ has Lipschitz constant $1$. If $K$ is a Cantor set then
$Im$ is surjective.

\item[(v)] If $X$ is connected, then $2^X$ is connected and if $K$ is a Cantor set then $\C(K,X)$ is connected. If $X$ is a
Cantor set then $2^X$ is a Cantor set and $\C(K,X)$, too, is totally disconnected.

\item[(vi)] For $h : K_1 \to K_2$ a continuous map of compact metrizable spaces, the map $h^* : \C(K_2,X) \to \C(K_1,X)$ defined by
$h^*(u) = u \circ h$ has Lipshitz constant at most $1$.  If $h$ is surjective then $h^*$ is an isometric inclusion and
$Im(h^*(u)) = Im(u)$.

\item[(vii)] For $f : X_1 \to X_2$ a continuous map of metric spaces, the maps $f_* : \C(K,X_1) \to \C(K,X_2)$ and
$f_* : 2^{X_1} \to 2^{X_2}$ are defined by $u \mapsto f \circ u$ and $A \mapsto f(A)$, respectively. If $f$ is
uniformly continuous or continuous, then the maps $f_*$ are uniformly continuous or continuous, respectively and
$f_*(Im(u)) = Im(f_*(u))$.

\item[(viii)] The following subsets  are closed.
\begin{equation}\label{1.4}
\begin{split}
INT \quad = \quad \{(A,B) \in 2^X \times 2^X: A \cap B \not= \emptyset \},  \hspace{2cm}\\
INC \quad = \quad \{(A,B) \in 2^X \times 2^X : A \subset B \}, \hspace{2.5cm} \\
EPS \quad = \quad \{(x,A) \in X \times 2^X : x \in A \}, \hspace{2.5cm}
\end{split}
\end{equation}
If $Y$ is a closed subset of $X$ then $2^Y$ is closed, regarded as a subset of $2^X$.
If $U$ is an open subset of $X$ then $2^U$ is open, regarded as a subset of $2^X$, and
$2^X \setminus 2^{X \setminus U} \ = \ \{ A : A \cap U \not= \emptyset \}$ is open.
\end{itemize}
\end{prop}

\proof The isometry and Lipschitz constant results in (i), (ii), (iv) and (vi) are easy to check. If $F \in 2^{2^X}$ and
$\{ x_n \}$ is a sequence in $\bigcup F$, then there exists $A_n \in F$ such that $x \in A_n$.
Since $F$ is compact we can assume by going to a subsequence that $\{ A_n \}$ converges to $A \in F$. By Lemma \ref{1.2}
 $\bigcup F \supset \bigcup_n \{ A_n \}  \cup A =  \overline{\bigcup_n \{ A_n \}}$ is compact and so $\{ x_n \}$
has a subsequence which converges to a point of $\bigcup F$. Hence, the latter is compact. When $K$ is a Cantor set, $Im$ is
surjective by Proposition \ref{prop1.1}. Since $\vee \circ i_{2^X} = 1_{2^X}$, $\vee$ is surjective. If $A$
contains $n+1$ distinct points at least $2 \ep$ apart then $A$ has distance at least $\ep$ from every set of cardinality
at most $n$, i.e. from the image of $i_{X,n}$.  Hence, the latter image is closed.

The density results in (iii) are clear. It follows that if $D$ is a countable dense set in $X$ then $FIN(D)$ is a countable dense
set in $2^X$ and $\bigcup_n \{ \C(\A_n,D) \}$ is a countable dense subset of $\C(K,X)$ when $K$ is a Cantor set and $\{\A_n\}$ is a
sequence of clopen partitions with mesh tending to zero. To approximate $u \in \C(K,X)$ we
apply uniform continuity of $u$ with respect
to any metric on $K$. If $K$ is not a Cantor set we apply Proposition \ref{prop1.1} to get
$h :K_1 \to K$ be a continuous surjection
with $K_1$ a Cantor set. Then $h^*$ is an isometric inclusion of $\C(K,X)$ onto a subset of the separable space $\C(K_1,X)$. Any
subset of a separable metric space is separable. If $X$ has no isolated points then neither does $FIN(X)$ or $\C(\A,X)$ for any
clopen partition $\A$. If $X$ is connected then each $\C(\A,X)$ is homeomorphic to $X^{(n)}$ where $n$ is the cardinality of $\A$ and
so is connected. By choosing a sequence $\{ \A_n \}$ of clopen partitions where each refines its predecessor and with the mesh
tending to zero we obtain a dense subset of $\C(K,X)$ as an increasing union of connected subsets. It follows that $\C(K,X)$ is
connected and so $2^X$ its image under $Im$ is connected as well.

If $X$ is a Cantor set then we can choose as metric on $X$ an \emph{ultra-metric} which satisfies the strengthening of the
triangle inequality:
$$d(x,y) \leq max (d(x,z),d(z,x)).$$ This is equivalent to saying that  $V_{\ep} = \{ (x,y) : d(x,y) < \ep \} $ is an open
equivalence relation for all $\ep > 0$. This implies that the open balls are clopen.  The ultrametric condition
 then holds for the induced metrics on $2^X$ and $\C(K,X)$ and so
each of these totally disconnected.  Since $2^X$ is compact and perfect it is a Cantor set, completing the proof of (v).

In (vii) the results are easy when $f$ is uniformly continuous. In general, if $\{ A_n \}$ is a sequence in
$2^{X_1}$ converging to $A$ then Lemma \ref{lem1.2}(a) implies that $A \subset \overline{\bigcup_n \{ A_n \}}$ and the
latter is  compact. So we can
use uniform continuity of the restriction of $f$ to this set to show that $\{ f_*(A_n) \}$ converges to $f_*(A)$.
Similarly, if $\{ u_n \}$ converges to $u$ in $\C(K,X_1)$ then $\{Im(u_n) \}$ converges to $Im(u)$ in $2^{X_1}$ and so
$u(K) \subset \overline{ \bigcup_n \{ u_n(K) \}}$ and the latter is compact. Again apply uniform continuity of $f$ on the
subset to show that $\{ f_*(u_n) \}$ converges to $f_*(u)$.

For (viii) it is easiest to proceed by checking that $A \cap B = \emptyset$, $A \cap (X \setminus B) \not= \emptyset$ and
$x \not\in A$ are open conditions.

$\Box$ \vspace{.5cm}

For us a dynamical system $(X,f)$ consists of a continuous map $f$ on a metric space $X$, i.e. $f : X \to X$. We are essentially interested in the dynamics of the  system $(2^X,f_*)$.

Such a study was first undertaken by Bauer and Sigmund \cite{bs}. The dynamics on $2^X$ is richer than the dynamics on $X$, which can be easily assessed from the results and examples in Bauer and Sigmund \cite{bs}, Glasner and Weiss \cite{gw}, John Banks \cite{ba} and Sharma and Nagar \cite{sn, sn1}. This provides us a strong motivation for the study of dynamics of $(2^X,f_*)$.

 In this article, we look into some known results - improving them,  some new concepts arising thereof, and answer some natural questions  with a view to further develop this study.

\section{Spaces of Continuous Maps, Spaces of Compact Subsets and Induced Dynamics}

For a system $(X,f)$ we call $(2^X,f_*)$ and $(\C(K,X),f_*)$ the \emph{induced systems}. The map $Im : \C(K,X) \to 2^X$ defines an
action map between the induced systems.  When $K$ is a Cantor set, it is a factor map.

 \begin{cor}\label{cor1.3a} Let $K$ be  Cantor set. If in a dynamical system $(X,f)$ the periodic
 points are dense then the periodic points
 are dense in the induced systems $(2^X,f_*)$ and $(\C(K,X),f_*)$. \end{cor}

 \proof Let $PER$ be the dense set of periodic points in $X$.  Any finite subset subset of $PER$ is a periodic point
 of $2^X$ and if $u \in C(\A,PER)$ with $\A$ a clopen partition of $X$, the $u$ is a periodic point of $\C(K,X)$.
 By (iv) above, $FIN(PER)$ is dense in $2^X$ and $\bigcup_{\A} \ \{ C(\A,PER) \}$ is dense in $\C(K,X)$.

$\Box$ \vspace{.5cm}

{\bfseries Remark:} Each periodic orbit is a fixed point in $2^X$.  So if $X$ is infinite and the periodic
points are dense, or more generally, if there are infinitely many periodic points in $X$ then there are infinitely
many fixed points in $2^X$.  It follows that if $X$ is compact and contains infinitely many periodic points
then the set of fixed points for the induced system on $2^X$ has an accumulation point.  Such a system cannot, for example, be
expansive.
\vspace{.3cm}

If $A, B \subset X$ we denote by $N(A,B)$  the \emph{hitting time set} $\{ n \in \N : A \cap f^{-n}(B) \not= \emptyset \}$.
Clearly, if $f_1$ and $ f_2$ are continuous functions on $X_1$ and $X_2$ respectively then for
 $A_1, B_1 \subset X_1, A_2, B_2 \subset X_2$  with respect to
 the product map $f_1 \times f_2$ on $X_1 \times X_2$ the hitting time set $N(A_1 \times A_2, B_1 \times B_2)$ is exactly
 $N(A_1,B_1) \cap N(A_2,B_2)$. The system is called \emph{topologically transitive} when $N(U,V) \not= \emptyset$ for every
 pair of  opene sets $U,V \subset X$. While $f$ need not be surjective, its image is clearly dense and
 so $U$ opene implies $f^{-n}(U)$ is opene for every $n \in \N$.

 The system is called \emph{weak mixing} when the product system $(X \times X, f \times f)$ is
 topologically transitive. It then follows that the $n$-fold product system $(X^{(n)},f^{(n)})$ is weak mixing for
 $n = 1,2,...$.  This follows by induction using the beautiful \emph{Furstenberg Intersection Lemma} \cite{fu}. We prove a
 strengthening due to Karl Petersen \cite{pk}.

 \begin{theo}\label{theo1.4} For a dynamical system $(X,f)$, assume that $N(U,V) \cap N(U,U) \neq  \emptyset$
 for every pair of opene sets $U,V \subset X$ and so the system is topologically transitive.
 For all opene sets $U_1,V_1,U_2,V_2 \subset X$ there exist opene sets
 $U_3,V_3 \subset X$ such that
 \begin{equation}\label{1.5}
 N(U_3,V_3) \quad \subset \quad  N(U_1,V_1) \cap  N(U_2,V_2). \hspace{2cm}
 \end{equation}
 In particular, the system is weak mixing.\end{theo}

 \proof $ N(U_1,V_1) \not= \emptyset$ implies there exists $n_1 \in\N$ such that $U_0 = U_1 \cap f^{-n_1}(V_1)$ is opene.
 $ N(U_0,U_2) \not= \emptyset$ implies there exists $n_2 \in\N$ such that $U = U_1 \cap f^{-n_1}(V_1) \cap f^{-n_2}(U_2)$ is opene.
Since $f$ is transitive, $f^{-n_1-n_2}(V_2)$ is opene.
\begin{equation}\label{1.6}
\begin{split}
N(U,U) \cap N(U,f^{-n_1-n_2}(V_2)) \quad \subset \hspace{2cm} \\
N(U_1,f^{-n_2}(U_2))\cap N(f^{-n_1}(V_1),f^{-n_1}(f^{-n_2}(V_2))) \hspace{.5cm}\\
= \quad  N(U_1,f^{-n_2}(U_2))\cap N(f^{n_1}(f^{-n_1}(V_1)),f^{-n_2}(V_2)) \\
\subset
\quad  N(U_1,f^{-n_2}(U_2))\cap N(V_1,f^{-n_2}(V_2)).\hspace{1.5cm}
\end{split}
\end{equation}
Fix $n_0 \in N(U_1,f^{-n_2}(U_2))\cap N(V_1,f^{-n_2}(V_2))$. With $n = n_0 + n_2$ the sets
$U_3 = U_1 \cap f^{-n}(U_2), V_3 = V_1 \cap f^{-n}(V_2)$ are opene.

Let $k \in N(U_3, V_3)$. Then $f^{-k}(V_3) \cap U_3 \neq \emptyset$. That is $f^{-k}(V_1) \cap f^{-n-k}(V_2) \cap U_1 \cap f^{-n}(U_2) \neq \emptyset$. Hence $k \in N(U_1, V_1) \cap N(f^{-n}(U_2), f^{-n}(V_2)) = N(U_1, V_1) \cap N(U_2, V_2)$.

As before, $N(f^{-n}(U_2),f^{-n}(V_2)) \subset N(U_2,V_2)$
and so (\ref{1.5}) follows.

$\Box$ \vspace{.3cm}

As a consequence, we obtain the following.

\begin{lem}\label{lem1.4a} If $f_i$ is a continuous map on a metric space $X_i$ for $i = 1,.., k$ and
the product system $(2^{X_1} \times ... \times  2^{X_k},(f_1)_* \times....\times (f_k)_*)$ is transitive then
the product system $(X_1 \times ... \times X_k, f_1 \times ... \times f_k)$ is weak mixing. \end{lem}

\proof: Given opene $U,V \subset X_1 \times ... \times X_n$ we show that $N(U,U)\cap N(U,V) \not= \emptyset$ and then apply
Theorem \ref{theo1.4}.

Choose opene $U_i, V_i \subset X_i$ such that $U_1 \times ... \times U_k \subset U$ and $V_1 \times ... V_k \subset V$.
Let $\bar U_i = 2^{U_i}$ and
$\bar V_i = \{ A \in 2^{X_i} : A \cap U_i \not= \emptyset \} \cap \{ A \in 2^{X_i} : A \cap V_i \not= \emptyset \}$.
Since $(f_1)_* \times ... \times (f_k)_*$ is transitive there exist $A_i \in \bar U_i$ and $n \in \N$
such that $((f_i)_*)^n(A_i) \in \bar V_i$ for $i = 1,...,k$. So there exist
$x_i, y_i \in A_i$ with $(f_i)^n(x_i) \in U_i, f^n(y_i) \in V_i$. Since $A_i \in 2^{U_i}$, $x_i, y_i \in U_i$.
Hence,
$n \in N(U_1 \times ... \times U_k,U_1 \times ... \times U_k) \cap N(U_1 \times ... \times U_k, V_1 \times ... \times V_k)$ and so
$n \in N(U,U) \cap N(U,V)$.

$\Box$ \vspace{.5cm}

An \emph{action map} $\pi : (X_1,f_1) \to (X_2,f_2)$ is a continuous map
$\pi : X_1 \to X_2$ such that $f_2 \circ \pi = \pi \circ f_1$.
When $\pi$ is surjective we call it a \emph{factor map} and say that $(X_2,f_2)$
is a factor of $(X_1,f_1)$.  It is easy to see that
a factor of a topologically transitive system is topologically transitive and
so a factor of a weak mixing system is weak mixing.

We now strengthen the result in \cite{ba, sn1}.

\begin{theo}\label{theo1.5} Let $K$ be a Cantor set. If $f$ is a continuous map on a metric space $X$
 then the following are equivalent:
\begin{itemize}
\item[(a)] $(X,f)$ is weak mixing.
\item[(b)] $(\C(K,X),f_*)$ is topologically transitive.
\item[(c)] $(2^X,f_*)$ is topologically transitive.
\end{itemize}
When these conditions hold, then the product system on $(X \times \C(K,X) \times 2^X,f \times f_* \times f_*)$ is weak mixing and
so each of the induced systems is weak mixing.
\end{theo}

\proof (a) $\Rightarrow$ (b): We prove that when $(X,f)$ is weak mixing, then $(X \times \C(K,X),f \times f_*)$
is topologically transitive. Notice that if $\A$ is a clopen partition of $K$ then $\C(\A,X)$ is an invariant set
and the subsystem $(\C(\A,X),f_*)$ is clearly isomorphic to the product system $(X^{(n)},f^{(n)})$
with $n$ the cardinality of $\A$. Hence, $(X \times \C(\A,X),f \times f_*)$ is topologically transitive for any clopen
partition $\A$. Now let $U_1,V_1 \subset X, U_2, V_2 \subset \C(K,X)$ be opene. By choosing the mesh of $\A$ fine enough
we obtain a clopen partition such that $U_2 \cap \C(\A,X), V_2 \cap \C(\A,X)$ are opene relative to $\C(\A,X)$.
Hence, $N(U_1 \times (U_2 \cap \C(\A,X)), V_1 \times (V_2 \cap \C(\A,X))) \subset N(U_1 \times U_2 , V_1 \times V_2)$ is nonempty
as required.  It then follows that the factor $(\C(K,X),f_*)$ is topologically transitive, proving (b).

(b) $\Rightarrow$ (c): Because $K$ is a Cantor set, $Im$ is a factor map and so transitivity on $2^X$ follows from
transitivity on $\C(K,X)$.

(c) $\Rightarrow$ (a): This is Lemma \ref{lem1.4a} with $k = 1$.

Above we showed that $(X,f)$ weak mixing implies $(X \times \C(K,X),f \times f_*)$ is topologically transitive.
Since the product system $(X \times X, f \times f)$ is weak mixing,
we see that the product system on $X \times X \times \C(K,X \times X)$ is
transitive. As $(\C(K,X \times X),(f \times f)_*)$ is isomorphic to $(\C(K,X) \times \C(K,X),f_* \times f_*)$ it follows that
the system on $X \times \C(K,X)$ is weak mixing and so the system on $X \times \C(K,X) \times \C(K,X)$ is weak mixing.  Applying
$Im$ on the third factor we see that $(X \times \C(K,X) \times 2^X,f \times f_* \times f_*)$ is a factor of the latter and
so is weak mixing as well.

$\Box$ \vspace{.3cm}

\begin{cor}\label{cor1.6}  If $(X,f)$ is weak mixing then the restriction of the product system to each of the closed,
invariant subsets $INT, INC \subset 2^X \times 2^X, EPS \subset X \times 2^X$ are weak mixing systems. \end{cor}

\proof Let $K$ be a Cantor set.  We will show that each of these systems is a factor of the weak mixing system
$(\C(K,X),f_*)$.

Let $\{K_0,K_1,K_2 \}$ be a partition of $K$ by three clopen sets. The map $\C(K,X) \to INT$ by
$u \mapsto (u(K_0 \cup K_1), u(K_0 \cup K_2))$ is clearly an action map. If $(A,B) \in INT$ we can choose continuous surjections
$u_0 : K_0 \to A \cap B, u_1 : K_1 \to A, u_2 : K_2 \to B$ and concatenate to obtain $u \in \C(K,X)$ mapping
to $(A,B)$.  The map $\C(K,X) \to INC$ by
$u \mapsto (u(K_0),u(K))$ is an action map.
If $(A,B) \in INC$ then choose surjections $u_0 : K_0 \to A, u_{12} : K_1 \cup K_2 \to B$ and concatenate.

Finally, fix $e \in K$ and map $\C(K,X)$ to $EPS$ by $u \mapsto (u(e),u(K))$. Again this is an action map.
If $(x,A) \in EPS$ then we can choose a continuous surjection $ u_0 : K \to A$. Let $y \in K$ with $u_0(y) = x$.
Since $K$ is homogeneous (e.g. there are Cantor sets which are topological groups) we can choose a homeomorphism
$h$ on $K$ such that $h(e) = y$. Then $u = u_0 \circ h$ maps to $(x,A)$.

$\Box$ \vspace{.5cm}

\section{Transitive Points for Induced Systems}

For the system $(X,f)$  and $x \in X$, the \emph{omega limit set of $x$}, denoted $\o f(x)$ is
the set of limit points of the orbit sequence.  Thus, $\o f(x) = Lim sup \{ \{ f^n(x) \} : n \in \N \}$.
We call $x$  a \emph{transitive point} when $\o f(x) = X$. If $\B$ is a base for the topology then
the set of transitive points, $Trans_f$ is  $\bigcap \{ \bigcup_{k \geq n} f^{-k}(U) : k \in \N, U \in \B \}$.
Clearly, if $X$ contains a transitive point then
the system is topologically transitive. In fact, every $N(U,V)$ is infinite. Furthermore, the space $X$ is separable and,
 unless $X$ consists of a single periodic orbit, the space
has no isolated points. Conversely, if $X$ is complete as well as separable
then topological transitivity
implies that $Trans_f$ is a dense $G_{\d}$ set. We apply the Baire category theorem to the above description of $Trans_f$ with
$\B$ a countable base.

Two points $x_1, x_2 \in X$ are called \emph{asymptotic} when
%eja display
$$Lim_{n \to \infty} \ d(f^n(x_1),f^n(x_2)) = 0.$$
  In that case,
$\o f(x_1) = \o f(x_2)$.  In particular, if $x_1 \in Trans_f$ then $x_2 \in Trans_f$.

  From now on we will assume that
%$K$ is a Cantor set and
$X$ is a complete,
separable,  metric space with no isolated points and so is a perfect Polish space. By Lemma \ref{lem1.2} and Proposition \ref{prop1.3} the same is
true of $2^X$ and $\C(K,X)$. Because $X$ is perfect, a point $x \in X$
is a transitive point if and only if the orbit sequence $\{ f^n(x) : n \in \N \}$ is dense in $X$.

A single periodic orbit with more than one point is not weak mixing. If $(X,f)$ is weak mixing then by Theorem \ref{theo1.5} the induced systems are weak mixing as well and so admit transitive points.
 The goal of this section is to examine the characteristics of these transitive points.

 \begin{df}\label{def2.1} Let $(X,f)$ be a dynamical system. A \emph{Kronecker subset} $L$ is a
 Cantor set contained in $X$ such that
 $\{ f^{n}|L : n \in \N \}$ is dense in $C(L,X)$. \end{df}
 \vspace{.5cm}

 This view of Kronecker subsets comes from Katznelson \cite{k}.

\begin{prop}\label{prop2.2} If $L$ is a Kronecker set for $(X,f)$ and $L_0$ is a nonempty clopen subset of $L_0$ then $L_0$ is a
Kronecker set for $(X,f)$. \end{prop}

\proof If $u_0$ is an arbitrary element of $C(L_0,X)$ and $\ep > 0$ then extend $u_0$ arbitrarily on the clopen
set $L \setminus L_0$ to obtain $u \in C(L,X)$.  There exists $n \in \N$ such that $d(u(x),f^n(x)) < \ep$ for all $x \in L$ and
so, a fortiori,  $d(u_0(x),f^n(x)) < \ep$ for all $x \in L_0$.

$\Box$ \vspace{.5cm}

\begin{theo}\label{theo2.3} $L$ is a Kronecker set for a system $(X,f)$ if and only if the inclusion map $u_L : L \to X$ is a transitive
point for the induced system $(\C(L,X),f_*)$. Equivalently,  $u$ is a transitive point for $(\C(K,X),f_*)$ with $K$ a Cantor set, if and only if
$u$ is injective with $Im(u) = u(K)$ a Kronecker set for $(X,f)$. \end{theo}

\proof $u$ is a transitive point for $(\C(K,X),f_*)$ if and only if $\{ f^n \circ u : n \in \N \}$ is dense in $\C(K,X)$. Since $\C(K,X)$
distinguishes points of $K$, it follows that $u$ must be an injective map in that case and so, by compactness, it is a
homeomorphism onto its image $L = u(K)$. So $ u^* : \C(L,X) \to \C(K,X)$ is an isomorphism between the induced systems. So
$u$ is a transitive point of $\C(K,X)$ if and only if $u_L$ is a transitive point of $\C(L,X)$. The $f_*$ orbit of $u_L$ is
$\{ f^{n}|L : n \in \N \}$ and so $u_L$ is a transitive point if and only if $L$ is a Kronecker subset.

$\Box$ \vspace{.5cm}

\begin{cor}\label{cor2.4} If $(X,f)$ is a dynamical system with $X$ a complete, separable, perfect metric space, then $(X,f)$
admits a Kronecker set if and only if the system is weak mixing. \end{cor}

\proof By Theorem \ref{theo2.3} $(X,f)$ admits a Kronecker subset if and only if the induced system on $\C(K,X)$ is topologically transitive.
By Theorem \ref{theo1.5} this occurs if and only if $(X,f)$ is weak mixing.

$\Box$ \vspace{.5cm}

In extending these results to the $n$-fold product, we recall the natural set isomorphisms
$(A^B)^C \cong A^{B \times C} \cong (A^C)^B$. We thus obtain natural isometric isomorphisms:
\begin{equation}\label{2.1}
\C(K, X^{(n)}) \quad \cong \quad \C(K \times \{ 1,...,n \},X) \quad \cong \quad \C(K,X)^{(n)}.
\end{equation}
Note that if $K$ is a Cantor set then $K \times \{ 1,..., n \}$ is a Cantor set with clopen partition
$\{ K \times \{ i \} : i = 1,...,n \}$. We will write $(u_1,...,u_n)$ for an element of any of these.

\begin{cor}\label{cor2.4a} Let $L \subset X^{(n)}$ with coordinate factors $L_i = \pi_i( L)$.

$L$ is a Kronecker set for $(X^{(n)},f^{(n)})$ if and only if the following three conditions holds
\begin{itemize}
\item The restriction of the coordinate projection $\pi_i :  L \to L_i$ is injective and so is a
homomorphism for $i = 1,...,n$.

\item The subsets  $\{ L_1,...,L_n \}$ of $X$ are pairwise disjoint.

\item The union $L_1 \cup ... \cup L_n$ is a Kronecker set for $(X,f)$.
\end{itemize}

Conversely, if $\{ L_1,..., L_n \}$ is a partition of a Kronecker set for $(X,f)$ by nonempty clopen subsets
then for any choice of homeomorphisms $u_i :K \to L_i$ of a fixed Cantor set $K$, the image $(u_1,...,u_n)(K) \subset X^{(n)}$
is a Kronecker set for $(X^{(n)},f^{(n)})$. \end{cor}

\proof  $L$ is a Kronecker set for $(X^{(n)},f^{(n)})$ if and only if it is the image of some transitive point $(u_1,...,u_n) \in \C(K,X^{(n)})$
for some Cantor set $K$. Via the isomorphisms of (\ref{2.1}) this is true if and only if $(u_1,...,u_n)$ is a transitive point
of $\C(K \times \{1,...,n\},X)$ which has image $L_1 \cup ... \cup L_n$ and is partitioned by $\{ L_1,..., L_n \}$. By
Theorem \ref{theo2.3} a transitive point for either induced system is an injective map and so is a homeomorphism onto its image.

$\Box$ \vspace{.5cm}

Now we consider $(2^X,f_*)$.

\begin{theo}\label{theo2.5}  Let $C$ be a transitive
point for $(2^X,f_*)$ with $X$ a complete, separable, perfect metric space.

\begin{itemize}
\item[(a)] If $g$ is a continuous function which commutes with $f$ and $C \cap g(C) \not= \emptyset$ then $g$ is the
identity map. In fact, if there exist $c_1, c_2 \in C$ with $c_2$ asymptotic to $g(c_1)$ then $g$ is the identity map.
The sets of the bi-infinite sequence $\{ f^{n}(C) : n \in \Z \}$ are pairwise disjoint. In particular,
$\{C, f(C), f^2(C),... \}$ is a pairwise disjoint sequence of transitive points for $(2^X,f_*)$.

\item[(b)] If $\mu$ is a invariant probability measure for $f$, then $\mu(C) = 0$.

\item[(c)] $C$ is nowhere dense.

\item[(d)] $C$ has infinitely many components. In particular, $C$ is infinite.

\item[(e)] Every point of $C$ is a transitive point for $f$.

\item[(f)] $C$ is a uniformly proximal set.  If  $x \in X$ then there is a sequence of iterates
$f^{j_i}(C)$ which converges to $\{ x \}$.
\end{itemize}
\end{theo}

\proof (a) Suppose $c_1, c_2 \in C$ with $g(c_1)$ asymptotic to $c_2$.  If $x \in X$ then there
exists a sequence $j_n$ such that $ f_*^{j_n}(C)$ converges
to $\{x \}$. Then $\{ f^{j_n}(c_1) \}$ converges to $x$ as does $\{ f^{j_n}(g(c_1)) = g(f^{j_n}(c_1)) \}$ because it is
asymptotic to $\{ f^{j_n}(c_2) \}$. Hence, $x = g(x)$. Since
$x$ was arbitrary, $g = 1_X$.

If two elements of the sequence $\{ f^{n}(C) : n \in \Z \}$ intersect then there exists $x \in C$ and a positive integer $k$
such that $f^k(x) \in C$. Hence, $f^k = 1_X$. This is impossible if $X$ is infinite and $f$ is transitive. Finally, all of
the points of the $f_*$ orbit of the transitive point $C$ are transitive points.

(b) Since the measure is invariant the sequence of disjoint sets $\{ f^{-k}(C) \}$ all have  measure $\mu(C)$. Since the sum is
finite $\mu(C) = 0$.

(c) We note that $x \in X$ is a  \emph{nonwandering point} if for every neighborhood $U$ of $x$, $f^n(U) \cap U \neq \emptyset$ for $n \geq 1$. The set of all nonwandering points is called the \emph{nonwandering set} of f. Any opene set is nonwandering and so the result follows from (a). That is, if $U, V \subset C$ are disjoint opene subsets
 then there exists a positive integer in
$N(U,V) \subset N(C,C)$ and this would contradict (a).

 %Alternatively, suppose $V_{\ep}(x) \subset C$.
% There exists $j > 0$ such that $f^j$ is $\ep$ close to $\{x \}$ and so
%$f^j(\overline{V_{\ep}(x)}) \subset V_{\ep}(x)$.  This means that $ \overline{V_{\ep}(x)}$ is an inward set for $f^j$ and
%so contains a nonempty attractor for $f^j$.  Since a weak mixing map is totally transitive, $f^j$ is transitive and so $X$ contains
%no $f^j$ attractors except for $X$ and $\emptyset$. It follows that $C$ has empty interior.

(d) Let $\{ x_1,.. , x_n \}$ be a set of $n$ distinct points and $2\ep < \min \ d(x_i, x_k)$ for $i \not= k = 1,...,n$. There
exists $j > 0$ such that $f^j(C)$ is $\ep$ close to $\{ x_1,.. , x_n \}$.  Hence, $\{ f^{-j}(V_{\ep}(x_i)) : i = 1,...,n \}$ is
a partition of $C$ into $n$ opene sets and hence the number of components is at least $n$.

(e) For every $U \subset X$ opene, there exists $j > 0$ such that $(f_*)^j(C) \in C(U)$, i.e. $f^j(C) \subset U$. Hence,
\begin{equation}\label{1}
j \in \bigcap_{x \in C} N(\{x \}, U).
\end{equation}

(f) The point $\{ x \}$ is a limit point of the $f_*$ orbit of $C$.

$\Box$ \vspace{.5cm}

\begin{lem}\label{lem2.6} Assume $C \in 2^X$.

(a)   If $A \in 2^X$ such that
\begin{equation}\label{3}
Lim_{j \to \infty} \ \max \{ d(f^j(a),f^j(C)) : a \in A  \} \quad = \quad 0,
\end{equation}
then $C \cup A$ is asymptotic to $ C$ with respect to $f_*$. In particular, if $C$ is a transitive point for $f_*$ then so is
$C \cup A$.

(b) Let $F$ be a finite subset of $X$ with $F \cap C = \emptyset$. If $B \in 2^X$ is perfect and
$(f_*)^{j_i}(C \cup F) \to B$ then $(f_*)^{j_i}(C) \to B$. If $C \cup F$ is a transitive point for $f_*$ then so is
$C$. \end{lem}

\proof (a) $ \max \{ d(f^j(a),f^j(C)) : a \in A  \}$ is the Hausdorff distance from $(f_*)^j(C \cup A)$ to $(f_*)^j(C)$.

(b) It suffices to show that every convergent subsequence of $ f_*^{j_i}(C)$ has limit $B$. Assume that
 subsequence of $ f_*^{j_i}(C)$ converges to $B_1$.  By going to a further subsequence we can assume that
$ f_*^{j_i}(F)$ converges to $F_1$.  Since the finite sets with cardinality at most that of $F$ form a closed subset of $2^X$,
$F_1$ is finite. Continuity of the map $\cup$ implies that $B_1 \cup F_1 = B$. Since $B$ is perfect and $B_1$ is closed, $B_1 = B$.
If $C \cup F$ is a transitive point and $B$ is a Cantor set transitive point then $B$ is in $\omega  f_*(C \cup F)$ and so is
in $\omega  f_*(C)$. Hence, the latter contains $\omega  f_*(B) = 2^X$.

$\Box$ \vspace{.5cm}

\begin{theo}\label{theo2.7} Let $(C_1,...,C_n)$ be a transitive point for $(f_*)^{(n)}$ on $(2^X)^{(n)}$.
\begin{itemize}
\item[(a)] The $C_i$'s are pairwise disjoint and so form a partition of $C = \bigcup_i C_i$. If $c_i \in C_i$ for
$i = 1,...,n$ then $(c_1,...,c_n)$ is a transitive point for $f^{(n)}$.

\item[(b)]The union $C = \bigcup_i C_i$
is a transitive point for $f_*$ and it contains $c_1,...,c_n$ such that $(c_1,...,c_n)$ is a transitive point for
$f^{(n)}$. The set of transitive points $C$ for which such a decomposition exists is dense in the set of transitive points
for $f_*$.
\end{itemize}
\end{theo}

\proof (a) Let $\{x_1,...,x_n \}$ be distinct points with $2 \ep > \min \ d(x_i,x_k) $ for $i \not= k = 1,...,n$.
There exists $j$ such that $(f_*)^j(C_i)$ is $\ep$ close to $\{x_i\}$ for $i = 1,...,n$ and so $\{ f^j(C_i) :i = 1,...n \}$
are pairwise disjoint and hence $\{C_i : i = 1,..., n \}$ are pairwise disjoint. If $c_i \in C_i$ then $f^j(c_i)$ is
$\ep$ close to $x_i$ for $i = 1,..., n$. This proves that $(c_1,...,c_n)$ is a transitive point of $f^{(n)}$.

(b) The map $(2^X)^{(n)} \to 2^X$ by $(A_1,...,A_n) \mapsto \bigcup_i A_i$ is surjective and continuous
 and so is a factor map from $(f_*)^{(n)}$ to $f_*$. A surjective continuous map takes dense sets to dense sets and
a factor map takes transitive points to transitive points.

$\Box$ \vspace{.5cm}

Call   $C$ an \emph{$n-$decomposable transitive point} when there
exists $(C_1,...,C_n)$  a transitive point for $(f_*)^{(n)}$ on $(2^X)^{(n)}$ such that $C = \bigcup_i C_i$.

\begin{cor}\label{cor2.8} If $C$ is an $n-$decomposable transitive point then it contains at least $n$ distinct
accumulation points. \end{cor}

\proof  If $C$ is $n-$decomposable then it admits a clopen partition $\{C_1,...,C_n \}$ with each $C_i$ infinite and
so with each containing an accumulation point.

$\Box$ \vspace{.5cm}

We will need the following routine result:

\begin{lem}\label{lem2.9} Let $L$ be a totally disconnected compact metric space.
If $\{C_1,..., C_n\} $ are pairwise disjoint, closed nonempty subsets of $L$ then there
exists a clopen partition $\{L_1,..., L_n\}$ of $L$ such that $C_i \subset L_i$ for $i = 1,..., n$.
\end{lem}

\proof Inductively we can for $i = 1,..., n-1$ choose $L_i$ a clopen
subset of $L$ which contains $C_i$ and which is disjoint from $L_j$ for $j < i$ and from $C_j$ for $j > i$. Then let
$L_n = L \setminus (L_1 \cup ... \cup L_{n-1})$.  Thus, $\{L_1,..., L_n\}$ is the required clopen partition of $L$.

$\Box$ \vspace{.5cm}

 A subset $E$ of $X$ is called {\emph{independent}} if any finite sequence $(e_1, e_2, \ldots, e_n)$ of points in $E$ is a transitive point in the product system $(X^{(n)},f^{(n)})$. Independent sets were studied by Iwanik \cite{iw}, and he proved the existence of large independent sets for weakly mixing systems. A subset of $X$ is called \emph{scrambled} when every pair of points in it is proximal and no pair of distinct points is
asymptotic. Independent sets are also wandering sets in $X$,  and  form a scrambled set in $X$.

\begin{theo}\label{theo2.10} Let $(X,f)$ be weak mixing and $L$ be a Kronecker subset of $X$.

\begin{itemize}

\item[(a)] $L$  is a transitive point for $f_*$ which is $n-$decomposable for every $n$. In fact, if
 $\{L_1,..., L_n \}$ is a partition of $L$ by nonempty clopen subsets then $(L_1,..., L_n)$ is a transitive point
 for $(f_*)^{(n)}$. If
 $c_1,..., c_n$ are distinct points in $L$ then $(c_1,..., c_n)$ is a transitive point for $f^{(n)}$, i.e. $L$ is an independent set.

 \item[(b)] If $C$ is any closed infinite subset of $L$ of $X$ then $C$ is a transitive point for $f_*$.
 If $\{C_1,...,C_n\} $ are pairwise disjoint, closed infinite subsets of $L$ then $(C_1,...,C_n)$ is a transitive point
 for $(f_*)^{(n)}$.

 \item[(c)] If $C$ is a closed subset of $L$ which contains  $n$ accumulation points then it is
 an $n-$decomposable transitive point. If it contains exactly $n$ accumulation points then it is not  $(n+1)-$decomposable.
 \end{itemize}
\end{theo}

\proof (a) $L$ is a Kronecker subset if and only if the inclusion map is a transitive point for $f_*$ on $\C(L,X)$. Since the map
$Im : \C(L,X) \to 2^X$ is a factor map, it takes
  transitive points  to transitive points. Applied to the inclusion map of $L$ itself this implies that $L \in 2^X$ is
 a transitive point.

If $\{L_1,..., L_n\}$ is a partition of $L$ by nonempty clopen subsets then by Proposition \ref{prop2.2} each is a Kronecker set.
Let $(B_1,..., B_n) \in (2^X)^{(n)}$. Choose $g_i \in \C(L_i,X) $ with $g_i(L_i) = B_i$. Concatenate to define
$g \in \C(L,X)$. There exists $j >0$ such that $f^j|L$ is close to $g$ and so $f^j(L_i)$ is close to $B_i$ for
$i = 1,...n$.  Hence, $(L_1,..., L_n)$ is a transitive point for $( f_*)^{(n)}$.
Thus,  $L$ is $n-$decomposable. Alternatively, one can
apply Corollary \ref{cor2.4a}.

If $c_1,..., c_n$ are distinct points of $L$ then there exists a partition $L_1,..., L_n$ of $L$ with $c_i \in L_i$
for $i = 1,..., n$. Since $(L_1,.., L_n)$ is a transitive point for $(f_*)^{(n)}$, $(c_1,..., c_n)$ is a transitive point
for $f^{(n)}$ by Theorem \ref{theo2.7} (a).

(b) If $c_1,..., c_n$ are distinct points of $L_0$ then they are distinct points of $L$ and so by (a) $(c_1,..., c_n)$ is a transitive
point of $f^{(n)}$.

If $\{C_1,..., C_n\} $ are pairwise disjoint, closed infinite subsets of $L$, then by Lemma \ref{lem2.9} there
is a clopen partition $\{ L_1,..., L_n \}$ with $C_i \subset L_i$ for $i = 1,..., n$.  Let $(B_1,..., B_n) \in (2^X)^{(n)}$ and $\ep > 0$.
There exists $M > 0$ and points $x_{ij} \in B_i, i = 1,..., n, j = 1,..., M$ so that $\{x_{ij} : j = 1,..., M \}$ is $\ep/2$ dense
in $B_i$ for $i = 1,..., n$. Now choose distinct points $c_{ij} \in C_i, i = 1,..., n, j = 1,..., M$. This is where we need
that each $C_i$ is infinite. By Lemma \ref{lem2.9} again there exist clopen partitions $\{L_{ij} :j = 1,..., M \}$ of $L_i$
for $i = 1,.., n$ such that $c_{ij} \in L_{ij}$ for all $i,j$.  Define $g \in \C(L,X)$ by $g(x) = x_{ij}$ for $x \in L_{ij}$.
There exists $k > 0$ such that
$f^k|L$ is $\ep/2$ close to $g$ and so $f^k(C_i)$ is $\ep/2$ close to $\{x_{i1},..., x_{iM} \}$ and so is $\ep$ close to $B_i$
for $i = 1,..., n$. This shows that $(C_1,..., C_n)$ is a transitive point for $( f_*)^{(n)}$.

(c) Let $c_1,..., c_n$ be  distinct accumulation points of $C$. By Lemma \ref{lem2.9} there is a clopen partition
$\{ C_1,..., C_n \}$ of $C$ with $c_i \in C_i$.  Hence, each $C_i$ is infinite. By (b) $(C_1,..., C_n)$ is a
$(f_*)^{(n)}$ transitive point and so $C$ is $n-$decomposable. If there are only $n$ accumulation points
then $C$ is not $(n+1)-$decomposable by Corollary \ref{cor2.8}.

$\Box$ \vspace{.5cm}

By combining Corollary \ref{cor1.6} with Theorem \ref{theo2.3} we sharpen the former.

\begin{cor}\label{cor2.11}  If $L$ is a Kronecker subset of $X$ and $\{ L_0, L_1,L_2 \}$ is
a clopen partition of $L$, then $(L,L_0)$ is a transitive point for $INC$ and $(L_0 \cup L_1, L_0 \cup L_2)$ is a transitive
point for $INT$. If $e \in L$ then $(e,L)$ is a transitive point for EPS. \end{cor}

%\proof The map $C(L,X) \to INC$ by $g \mapsto (g(L),g(L_0))$ is onto:  Given $(A,B) \in INC$ first choose $g_0 :L_0 \to B$ onto
%and then choose $g_12 : L_1 \cup L_2 \to A$ onto.  Concatenate to define $g$. The transitive point $inc_L$ of $C(L,X)$ is mapped
%to $(L,L_0)$.  The map $C(L,X) \to INT$ by  $g \mapsto (g(L_0 \cup L_1),g(L_0\cup L_2))$ is onto:  Given $(A,B) \in INT$
% first choose $g_0 :L_0 \to A \cap B$ onto
%and then choose $g_1 : L_1  \to A$ onto and $g_2 : L_2 \to B$ onto. Concatenate to define $g$.
%The transitive point $inc_L$ of $C(L,X)$ is mapped
%to $(L_0 \cup L_1,L_0 \cup L_2)$.

$\Box$ \vspace{.5cm}

\begin{cor}\label{cor2.12}
(a) There exists a transitive point $C$ for $f_*$ which contains a single accumulation point $c^* \in C$,
and so $C \setminus \{ c^* \}$ is a set of isolated points in $X$. In particular, $C$ is countable. Such a set $C$ is
never $2-$decomposable.  However, it can be chosen so that if $c_1,..., c_n$ are distinct points of $C$ then
$(c_1,..., c_n)$ is a transitive point for $f^{(n)}$.

(b) If $A$ is any infinite, compact, totally disconnected metric space there exists a transitive point $C$ for $f_*$ which
is homeomorphic to $A$.
\end{cor}

\proof (a) Any set of isolated points in a separable metric space is countable. A set with a single accumulation point is not
$2-$decomposable by Corollary \ref{cor2.8}.
 %If $C \setminus \{c^*\}$ is a set of isolated
%isolated points and $C_1, C_2$ form a partition of $C$ be closed sets with $c^* \in C_1$ then $C_2$ is finite and so
%cannot be a transitive point for $ \bar f$. A fortiori, $(C_1,C_2)$ is not a transitive point for $\bar f \times \bar f$.

Let $L$ be a Kronecker subset and let $\{ c_k : k = 1,2,.. \}$ be a sequence of distinct points in $L$ converging
to a point $c^*$ in $L$. Let $C$ be the points of this sequence together with the limit point. Since the sequence
is convergent, the points of $C \setminus c^*$ are isolated. As it is an infinite subset of $L$ it is a transitive point
for $f_*$ by Theorem \ref{theo2.10} (b).

%Let $B \in 2^X$ and $\ep > 0$ choose $x_1,...,x_n$ points of $B$ so that $\{x_1,...,x_n \}$ is an $\ep/2$ dense subset of
%$B$. Choose $\{ V_i : i = 1,...,n-1 \}$ disjoint clopen subsets of $L$ such that $V_i \cap C = \{c_i \}$. Let $V_n = L \setminus
%(\bigcup_{i=1}^{n-1} V_i)$.  Define $g \in C(L,X)$ by $g(x) = x_i$ for $x \in V_i$. There exists $j > 0$ such that
%$f^j|L$ is $\ep/2$ close to $g$ and so $f^j(C)$ is $\ep/2$ close to $\{x_1,...,x_n \}$ and so is $\ep$ close to $B$.

(b) The product $A \times L$ is perfect as well as totally, disconnected and compact metric.  Hence, $A \times L$ is homeomorphic
to $L$.  Hence, $L$ contains a homeomorphic copy of $A$.  Such a set is an infinite closed subset of $L$ and so is a transitive point
by Theorem \ref{theo2.10} (b) again.

$\Box$ \vspace{.5cm}
%
%{\bfseries Remark:} If we take $C_1,...,C_n$ disjoint convergent sequences in $L$ and then take the union, we obtain
%a countable transitive point $C$ for $\bar f$ with $n$ accumulation points and which is $n-$decomposable but not
%$n+1-$decomposable.
%\vspace{.5cm}

A subset of $X$ is called  \emph{strongly scrambled} when every pair of points in it is proximal and recurrent. Any Kronecker set
is strongly scrambled. The map $f$ is
called \emph{uniformly rigid} when some sequence of iterates $f^{j_i}$ converges uniformly to the identity on $X$.
In that case, every pair of points is recurrent and so no pair of distinct points is asymptotic.  Glasner and Maon \cite{gm} have
constructed uniformly rigid, weak mixing, minimal homeomorphisms on the torus.  If $f$ is uniformly rigid then every transitive point
$C$ for $f$ is strongly scrambled.

We note that Hernandez , King  and Mendez-Lango  \cite{hkm} have constructed  Cantor sets with dense orbit in $(2^X,f_*)$.

\section{Examples}

We first look into an example of a Kronecker set.

Let  $X = \{0,1\}^\mathbb{N}$ along with the shift map $\sigma$ be the one-sided shift  space.

Consider the induced system $(2^X, \sigma_*)$.

Let $\mathcal{L}(X)$ denote the language of $X$, and we define $W_n$ to be the set of all words of length $n$ in $\mathcal{L}(X)$. So we have,
$$W_1 = \{w^{(1)}_1, w^{(1)}_2\}, \  W_2 = \{w^{(2)}_1, w^{(2)}_2, w^{(2)}_3, w^{(2)}_4\}, \ldots, W_n = \{w^{(n)}_1, \ldots, w^{(n)}_{{2^n}}\}$$

where $w^{(j)}_i$ can be considered lexicographically in each $W_j$, i.e. in each $W_j$, $w^{(j)}_1 < w^{(j)}_2 < \ldots < w^{(j)}_{2^j-1} < w^{(j)}_{2^j}$.

Also, we consider $S_n$ - the group of permutations on $n$ elements.

And let $\mathcal{U}^{n,k} \subset \mathcal{P}(W_n)$ for $1 < k \leq 2^n$, be the collection of all sets in the power set of $W_n$ with $k$ elements. We note that there will be ${_{2^n}}C_k$ of them, which we can enumerate as $U^{n,k}_1, U^{n,k}_2, \ldots, U^{n,k}_{{_{2^n}}C_k}$. We note that $|U^{n,k}_j| = k$, and for each $u^{n,k}_j \in U^{n,k}_j$, we have $|u^{n,k}_j| = n$

\vskip .5cm

Define $C \in 2^X$ as

\vskip .5cm

$ C = \{ w^{(1)}_{\rho_{1_1}(1)}w^{(1)}_{\rho_{1_1}(2)}w^{(1)}_{\rho_{1_2}(1)}w^{(1)}_{\rho_{1_2}(2)}$  $w^{(2)}_1 w^{(2)}_2 \ldots w^{(2)}_4$
  $ u^{2,2}_1 u^{2,2}_2 \ldots u^{2,2}_6$ $u^{2,3}_1 u^{2,3}_2 u^{2,3}_3 u^{2,3}_4$ $u^{2,4}_1$
  $w^{(3)}_{\rho_{3_1}(1)}w^{(3)}_{\rho_{3_1}(2)} \ldots w^{(3)}_{\rho_{3_1}(2^3)}  \ldots$
    $w^{(3)}_{\rho_{3_{2^3}}(1)}w^{(3)}_{\rho_{3_{2^3}}(2)} \ldots w^{(3)}_{\rho_{3_{2^3}}(2^3)}$

   $w^{(4)}_1   w^{(4)}_2    w^{(4)}_3  w^{(4)}_4  \ldots w^{(4)}_{2^4}$ $u^{4,2}_1 \ldots u^{4,2}_{{_{2^4}}{C_2}} \ldots u^{n,2^n}_1$
    $ w^{(5)}_{\rho_{5_1}(1)}w^{(5)}_{\rho_{5_1}(2)} \ldots w^{(5)}_{\rho_{5_1}(2^5)}$

    $\ldots w^{(5)}_{\rho_{5_{2^5}}(1)}w^{(5)}_{\rho_{5_{2^5}}(2)} \ldots w^{(5)}_{\rho_{5_{2^5}}(2^5)}$
     $ w^{(6)}_1 \ldots$
     $u^{2n,2}_1 \ldots u^{2n,2}_{{_{2^{2n}}}{C_2}} \ldots u^{2n,2^{2n}-1}_{{_{2^{2n}}}{C_{2^{2n}-1}}} u^{2n,2^{2n}}_1$

      $w^{(2n+1)}_{\rho_{{2n+1}_1}(1)} \ldots w^{(2n+1)}_{\rho_{{2n+1}_1}(2^{2n+1})} \ldots
       / u^{n,k}_j \in U^{n,k}_j$, $\rho_{k_j} \in S_{2^k}$ for each $n, k, j \in \mathbb{N} \}$,

\vskip .5cm

i.e. $C$ consists of all sequences with consecutive blocks of all length $1$ words arranged according to all permutations in $S_2$, followed by  consecutive blocks of all length $2$ words arranged lexicographically, followed by  words of length $2$ from each set containing $2, 3, 4$ of these words of length $2$, followed by consecutive blocks of all length $3$ words arranged according to all permutations in $S_8$, \ldots, followed by consecutive blocks of all length $2n$ words arranged lexicographically, followed by  words of length $2n$ from each set containing $2, 3, \ldots, 2^{2n}$ of these words of length $2n$, followed by consecutive blocks of all length $2n+1$ words arranged according to all permutations in $S_{2^{2n+1}}$, $\ldots$.

\vskip 1cm

We claim that $C$ is compact. It is enough to see that $C$ is closed.

Let $x = x_1 \ldots x_p \ldots \notin C$. This means that $x$ cannot be realized as a sequence with a consecutive arrangement of  blocks of all length $n$ words for $n \in \mathbb{N}$ as done in $C$. Let $t$ be the smallest index where this difference can be realized, and let $j$ be largest even integer such that

$r = t - ({2 \cdot 1 + 2^2 \cdot 2+ 2 \cdot \sum \limits_{k=2}^4 {_4}{C_k} + 2^3 \cdot 3 + \dots +2^j \cdot j + j \cdot \sum \limits_{k=2}^{2^j} {_{2^j}}{C_k}}) > 0$, and $t - r < 2^{j+1} \cdot j$.

Then, $x \in [x_1 \ldots x_h]$, but $C \cap [x_1 \ldots x_h] = \phi$, when $ h = 2 \cdot 1 + 2^2 \cdot 2+ 2 \cdot \sum \limits_{k=2}^4 {_4}{C_k} + 2^3 \cdot 3 + \dots +2^j \cdot j + j \cdot \sum \limits_{k=2}^{2^j} {_{2^j}}{C_k} + 2^{j+1} \cdot (j+1)$.

\vskip 1cm

We see that $C$ is transitive in $(2^X, \sigma_*)$.

Take any basic open set $< A_1, \ldots, A_q >$ in the Vietoris topology on $2^X$. We can always have a $p \in 2\mathbb{N}$ large enough so that the cylinders $[a^i_1 \ldots a^i_p] \subset A_i, \ 1 \leq i \leq q$. For some $j \in {1, \ldots, 2^p}$, we have $U^{p,q}_j = \{ a^1_1 \ldots a^1_p, \ldots, a^q_1 \ldots a^q_p \}$. Then we have $\sigma_*^h (x)$ $\subset$ $< [a^1_1 \ldots a^1_p], \ldots, [a^q_1 \ldots a^q_p]>$ $\subset$ $< A_1, \ldots, A_q >$, where $h =  2 \cdot 1 + 2^2 \cdot 2+ 2 \cdot \sum \limits_{k=2}^4 {_4}{C_k} + 2^3 \cdot 3 + \dots +2^p \cdot p + p \cdot \sum \limits_{k=2}^{j-1} {_{2^p}}{C_k}$.
 \vskip 1cm

We note that each $c \in C$ is a transitive point for $(X, \sigma)$, and for every $n \in \N$,  $ C = C_1 \cup \ldots \cup C_n$ such that $(C_1,..., C_n)$ is a $(\sigma_*)^{(n)}$ transitive point and so $C$ is $n-$decomposable.

$\Box$ \vskip 1cm

Now we specialize to $X = \{0,1\}^{\Z}$ with $ \sigma$ the shift homeomorphism.

\begin{theo} \label{theo3.1} There exists a transitive point $C$ for $\sigma_*$ which
is a countable set with a single accumulation point $c^* \in C$.
In addition, any pair of points $c_1, c_2 \in C$ is asymptotic.
In particular, $(c_1,c_2)$ is never a transitive point for $\sigma \times \sigma$.
\end{theo}

\proof Call a set of finite words an EL set (equal length) if they all have the same length.  Let $\{A_1, A_2, .... \}$
be a sequence which counts all those EL sets on the alphabet $\{0, 1\}$  whose common length is odd.  Let $2\ell_k + 1$
be the common length of the words in $A_k$. In each $A_k$ choose a particular word $a_k^*$ which we call the
\emph{special words}.  For example, using lexicographic
ordering we could choose the first word in the set $A_k$. Let $N_k = \Sigma_{j=1}^k (2 \ell_k + 1)$.

We now construct inductively a sequence of EL sets $\{ B_1,... \}$ for which $N_k$ is the common length of the words in
$B_k$ and there is a distinguished element $b^*_k \in B_k$. Each member of $B_k$ is a concatenation of blocks
of length $2 \ell_1 + 1,..., 2 \ell_k + 1$. We refer to these as the first block, the second, up to the $k^{th}$ block.

Let $B_1 = A_1$ and $b^*_1 = a^*_1$.

Inductively, we concatenate, defining $B_{k+1}$ to be the union of the two sets $\{ b^*_k \} \cdot A_{k+1}$ and
$B_k \cdot \{a^*_{k+1} \}$. The intersection of these two sets is the single word $b^*_k \cdot a^*_{k+1} $ which
we define to be $b^*_{k+1}$.

Thus, $b^*_{k+1} = a^*_1 \cdot a^*_2 \dot ... a^*_k \cdot a^*_{k+1}$.

Inductively it is clear that for each word in $B_{k+1}$ at most one block is not a special word.

Define $c \in C$ when $c_i = 0$ for all $i \leq 0$ and when $c_{[1,...,N_k]} \in B_k$ for all $k$. Let
$c^*$ be the special element of $C$ with $c^*_i = 0$ for all $i \leq 0$ and $c^*_{[1,...,N_k]} = b^*_k$.

If $a_k \in A_k \setminus \{ a^*_k \}$ and the $k^{th}$ block of $c \in C$ is $a_k$ then every other block on (the
positive side) is a special word. Thus, $c$ agrees with $c^*$ except in finitely many places and it is clear
 that $c$ is an isolated point of $C$. Thus, the points
of any finite subset of $C$ all agree with $c^*$ after some finite coordinate level.  It follows that any two points are asymptotic.
  This in turn implies that no pair in $C$ defines a transitive point for $\sigma \times \sigma$.

Finally, we observe that the words from the $-\ell_{k+1}$ to the $+\ell_{k+1}$ coordinate occurring in the
set $\sigma^{N_k + \ell_{k+1} +1}(C)$
are exactly the words of $A_{k+1}$ centered in the middle.  Since the sequence $\{A_k \}$ lists all EL sets of words of
odd length it follows that the iterates $\sigma_*^j(C)$ approach every element of $2^X$ arbitrarily closely.  That is,
$C$ is a transitive point for $\sigma_*$.

$\Box$ \vspace{.5cm}

We can obtain  Cantor set examples by using the \emph{asymptotic extension technique}.  For $B \in 2^X$ with $X = \{0,1 \}^{\Z}$
let $B_N = \pi_N^{-1}(\pi_N(B))$ where $\pi_N :  \{ 0,1 \}^{\Z} \to \{0,1 \}^{[N,\infty)}$ is the coordinate projection. It is
clear that $B$ and $B_N$ are asymptotic points of $2^X$, with respect to $\sigma_*$,  since they agree beyond the $N^{th}$ coordinate.
Furthermore, every point of $B_N$ is asymptotic to a point of $B$.  In particular, if $B$ is a transitive point for $ \sigma_*$ then
each $B_N$ is as well.  %If $B $ is $n-$decomposable then each $B_N$ is since $B = B_1 \cup ... \cup B_n$ implies
%$B_N = (B_1)_N \cup ... \cup (B_n)_N$.

\begin{cor} \label{cor3.2} There exists a transitive point $C$ for $\sigma_*$ which
is a Cantor set such that
 any pair of points $c_1, c_2 \in C$ is asymptotic.
In particular, $(c_1,c_2)$ is never a transitive point for $\sigma \times \sigma$.
\end{cor}

\proof Let $B$ be the countable set constructed for Theorem \ref{theo3.1} so that every point of $B$ is asymptotic to $c^*$.
Clearly, every point of $C = B_N$ is asymptotic to $c^*$. Since the coordinates
below $N$ are arbitrary $B_N$ is a Cantor set.

$\Box$ \vspace{.5cm}

{\bfseries Remark:} By extending the idea of Theorem \ref{theo3.1} one can construct $B$ to be a disjoint union of sets
$B_1, .., B_n$ each an infinite countable set contained
in a single asymptotic class and so that $(B_1,...,B_n)$ is a transitive point for $(\sigma)^{(n)}$. Then $B_N$ is a Cantor set
$n$-decomposable transitive point which meets only $n$ asymptotic classes and so is not $n+1-$decomposable.
\vspace{.5cm}

Using this asymptotic extension idea in a different context we can obtain examples of transitive points in $2^X$ with nontrivial
components.

We specialize to $X = \R^2/\Z^2$ with $\pi : \R^2 \to \R^2/\Z^2$ the canonical projection. We define the metric
$d$ on $\R^2/\Z^2$ by $d(a + \Z^2,b + \Z^2) = min \{ |a - b + z|: z \in \Z^2 \}$.  Hence, $\pi$ has Lipschitz constant
1.  Any closed  $C \subset \R^2/\Z^2$  of diameter less than $1$ is contained in a ball in $\R^2/\Z^2$ and by lifting the ball
we obtain a closed
 $C^+ \subset \R^2$ such that $\pi : C^+ \to C$ is a homeomorphism.  If $C \subset \R^2/\Z^2$ is a Cantor set of any diameter, we can
 choose a clopen partition of $C$ by sets of diameter less than one and lift each separately to obtain a Cantor set
 $C^+ \subset \R^2$ such that $\pi : C^+ \to C$ is a homeomorphism.

Now let $t$ be a Thom torus map on $\R^2/\Z^2$, e.g. the diffeomorphism induced by the linear map $T$
with matrix $\begin{pmatrix} 3 & 1 \\ 5 & 2 \end{pmatrix}$.
This is an Anosov diffeomorphism and so is a factor of the shift map.  In particular, it is weak mixing. The eigenvalues
$\l_{\pm} = (5 \pm \sqrt{21})/2$ satisfy $4 < \l_+ < 5$ and $0 < \l_- < \frac{1}{2}$. Let $v_+$ and $v_-$ be  corresponding
eigenvectors of unit length. The lines parallel to $v_-$ and to $v_+$ project to the stable and unstable foliations. In particular,
for any $a \in \R^2$ and any real number $s$
\begin{equation}\label{3.1}
|T^j(a) - T^j(a + s v_-)| \quad = \quad |s| (\l_-)^j \quad < \quad |s| 2^{-j}.
\end{equation}
In particular, all points on the stable manifold of $\pi(a)$, i.e. $\pi(a + \R v_-)$, are asymptotic with respect to $t$.

Let $L$ be a Kronecker subset of $\R^2/\Z^2$ for $t$. There exists a sequence $j_k \to \infty$ such that $t^{j_k}|L$ converges to
the inclusion map of $L$ into $\R^2/\Z^2$.  Hence, no two distinct points of $L$ are asymptotic. Hence, the stable manifolds
of any two distinct points in $L$ are disjoint.  Let $L^+$ be a Cantor set in $\R^2$ such that $\pi|L^+ \to L$ is a homeomorphism.
Define the map $Q : L^+ \times \R \to \R^2/\Z^2$ by $Q(a,s) = \pi(a + s v_-)$. The map Q is clearly continuous. It is injective
because the $v_-$ lines through distinct points of $L^+$ do not intersect.  For $M > 0$ let $C_M = Q(L^+ \times [-M,M])$.
By compactness, the restriction of $Q$ is a homeomorphism from $L^+ \times [-M,M]$ onto $C_M$. From equation (\ref{4})
it follows that the Hausdorff distance between $T^j(L^+)$ and $T^j(L^+ \ + \ [-M,M]v_-)$ converges to zero. Hence, the distance
between $ t_*^j(L)$ and $ t_*^j(C_M)$ converges to zero. That is $L$ and $C_M$ are asymptotic points of $2^{X}$.
Since $L$ is a transitive point, $C_M$ is as well.  Thus, we have proved:

\begin{theo} \label{theo3.3} There exists a transitive point $C$ for $t_*$ which is homeomorphic to $L \times [0,1]$ where
$L$ is a Cantor set.  In particular, every component of $C$ is an interval. \end{theo}

$\Box$ \vspace{.5cm}

{\bfseries Remark:}  Note that $\{ C_M : M = 1,2,... \}$ is an increasing sequence of transitive points for $t_*$ whose
union is dense in $\R^2/\Z^2$.
 \vspace{.5cm}

In examples of this sort any pair of points which lie in the same component of the transitive point $C$ are asymptotic. The question
arises whether this is always true.  We obtain a partial result suggesting that this is true.

For a compact metric space $X$ define the \emph{component diameter} to be
\begin{equation}\label{3.2}
Cdiam(X) \quad = \quad sup \{ d(x,y) : x  \ \mbox{and} \ y \ \mbox{lie in the same component of} \ X \}.
\end{equation}
Equivalently, $Cdiam(X)$ is the supremum of the diameters of the components of $X$.

\begin{lem}\label{lem3.4} If $\{ B_i \}$ is a sequence in $2^X$ converging to a totally disconnected $C \in 2^X$ then
\begin{equation}\label{6}
Lim_{i \to \infty} \ Cdiam(B_i) \quad = \quad 0. \hspace{1cm}
\end{equation}
\end{lem}

\proof  Given $\ep > 0$ choose a clopen partition of $C$ by elements of diameter less than $\ep$ and enlarge each to
an open subset of $X$.  That is, we can obtain a disjoint opene sets $U_1,...,U_n$ whose union $U$ contains $C$. There
exists $N$ so that $i \geq N$ implies $B_i \subset U$ and so each component is entirely contained in one of the $U_i$'s.
Hence, $i \geq N$ implies that $Cdiam(B_i) < \ep$.

$\Box$ \vspace{.5cm}

%eja included this conjecture
\emph{We conjecture that if $C$ is a transitive point for $f_*$ on $2^X$, then
$$Lim_{i \to \infty} \ Cdiam((f_*)^{i}(C)) \quad = \quad 0.$$}
The closest we can get is the following result

\begin{theo}\label{theo3.5} Assume $f$ is a homeomorphism. $C$ is a
transitive point for $f_*$ on $2^X$ and $B \in 2^X$ then there exists
a sequence $\{ j_i \}$ of integers tending to infinity such that:
\begin{equation}\label{7}
\begin{split}
Lim_{i \to \infty} \ \{ (f_*)^{j_i}(C) \} \quad = \quad  B, \qquad \mbox{and} \\
Lim_{i \to \infty} \ Cdiam((f_*)^{j_i}(C)) \quad = \quad 0. \hspace{1cm}
\end{split}
\end{equation}
\end{theo}

\proof  Let $L$ be a Cantor set transitive point for $f_*$ and choose sequences $n_p, m_q$ so that
$(f_*)^{n_p}(C) \to L$ and $(f_*)^{m_q}(L) \to B$ as $p, q \to \infty$.  Diagonalize.  That is, given $i > 0$ choose
$q_i$ so that $d(B,(f_*)^{m_{(q_i)}}(L)) \ < \ 1/2i$ and then $p_i$ so that
$d((f_*)^{m_{(q_i)}}(L),(f_*)^{m_{(q_i)} + n_{(p_i)}}(C)) \ < \ 1/2i$ and $Cdiam((f_*)^{m_{(q_i)} + n_{(p_i)}}(C)) < 1/i$.
Because $f$ is a homeomorphism $(f_*)^{m_{(q_i)}}(L)$ is a Cantor set and so Lemma \ref{lem3.4} implies that $p_i$ exists.
Let $j_i = m_{(q_i)} + n_{(p_i)}$.

$\Box$ \vspace{.5cm}

If $H(K)$ is the homeomorphism group of the Cantor set $K$ then $H(K)$ acts on $\C(K,X)$ by $tg = g \circ t^{-1}$. This action
commutes with $f_*$.
The map $Im : \C(K,X) \to 2^X$ is constant on the $H(K)$ orbits. Let $\C_{in}(K,X)$ be the set of injective continuous maps.
This is a dense $G_{\d}$ subset of $\C(K,X)$ and it is invariant under $f_*$. It maps to $CANTOR \subset 2^X$ the dense $G_{\delta}$
subset of Cantor sets in $X$. Furthermore the fibers of the map $\C_{in}(K,X) \to CANTOR$ are exactly the $H(K)$ orbits.

%eja I changed the title to get compactness into it
\section{Minimality and Distality for the Induced Dynamics of Compact Systems}

The system $(X,f)$
is \emph{minimal} if every $x \in X$ is a transitive point. We note that $X$ is always a fixed point for $(2^X,f_*)$, and so the system $(2^X,f_*)$ can never be minimal. Here we talk of minimality in a weaker sense. In this section $X$ will always be a compact metric space.

\begin{df}
The system $(X,f)$ is said to be \emph{backward minimal}
 if the negative semiorbit of $x$, $ O^-(x)= \{y \in X :
f^n(y)=x\ $ for some $n \in \N \}$, is dense in $X$ for every $x\in X$.

Equivalently, $(X,f)$ is  \emph{backward minimal} if
for every  opene $U$, $X = \bigcup \limits_{n=1}^\infty
f^n(U) $.
\end{df}

{\bfseries Remark:} We can immediately make the following conclusions:

1. A minimal system is always backward minimal and if $f$ is a homeomorphism, then minimal and backward minimal are equivalent.

2. A backward minimal system is always topologically transitive.

3. A system $(X,f)$ is called \emph{exact} if for every opene $U$ in $X$, there exists an $n \in \N$ such that $f^n(U) = X$. An exact system is always backward minimal.

A simple example of a backward minimal system is the one-sided shift space on finite symbols. It is easy to see that this shift space will be exact but not minimal. Also, the two-sided shift space on finite symbols, where the shift map is now a homeomorphism, fails to be backward minimal. The irrational rotation is a  minimal system so is backward minimal but is not exact.

\begin{theo} A backward minimal system $(X,f)$ need not be weakly mixing. In particular, in this case $(2^X,f_*)$ need not be topologically transitive.
\end{theo}

\proof Let $X$ be the closure  in $\{0,1\}^{\N}$ of the sequences of the form
$$0^k10^{3^{n_1}}10^{3^{n_2}}10^{3^{n_3}}1  \ldots 10^{3^{n_i}}1 \ldots$$
 where $k, n_i \in \N$. Consider the shift map $\sigma$ on $X$. We see that for any cylinder, a finite image will be a cylinder of some sequences that start with $1$, and consequently the images of this cylinder will cover the whole of $X$. So this system $(X, \sigma)$ is backward minimal.

However,  there is no  $n \in \N$ for which $\sigma^n[001] \cap [1] \neq \phi$ and also $\sigma^n[010] \cap [1] \neq \phi$ simultaneously. So our system $(X, \sigma)$ is not weakly mixing. By Theorem \ref{theo1.5}, $(2^X, \sigma_*)$ cannot be topologically transitive.

 $\Box$

\begin{theo}  For a dynamical system $(X,f)$ the following are equivalent:
\begin{itemize}
\item[(i)] $(X,f)$ is exact.
\item[(ii)] $(2^X,f_*)$ is exact.
\item[(iii)] $(2^X,f_*)$ is backward minimal.
\item[(iv)] The negative orbit $ O^-( \{ X \} ) \ $ of the point $X$ is dense in $2^X$.
\item[(v)] $i_X(X) \subset 2^X$ is contained in the closure  of $ \ \bigcup_{i \geq 0} (f_*)^{-i}( \{ X \} )$ in $2^X$.
\end{itemize}
\end{theo}

%eja moving the previous discussion to make a proof
\proof  It is obvious that (ii) $\Rightarrow $ (iii) $ \Rightarrow $ (iv) $\Rightarrow $ (v).

(v) $\Rightarrow$ (i): We are assuming that the backward orbit of the point $X \in 2^X$ meets every neighborhood
of  $\{ x \} \in 2^X$ for any  $x \in X$. Thus, any opene subset of $X$ contains a closed set $A$ which
 maps onto $X$ under some iterate of $f$. That is, $f$ is exact.

(i) $\Rightarrow$ (ii):   Since $f$ is exact, it is surjective and
 so $f^n(U) = X$ implies $f^i(U) = X$ for  $i \geq n$. Now, given $A \in 2^X$,
 take a Vietoris neighbourhood of $A$ i.e.  a list of closed sets $C_1,...,C_k$  covering $A$ and such that each interior meets
 $A$.  There exists $n \in \N$ so that $f^n(Int C_i) = X$ for $i = 1,...,k$. For any $B \in 2^X$, let $D_i = (f^{-n}|C_i)(B)$.
  Each $D_i \subset C_i$ and
 meets $Int C_i$. Hence,
 $A_1 = \bigcup \limits_{i=1}^k D_i$ is in the Vietoris neighbourhood  of $A$, while $f^n(A_1) = B$. Hence, $f_*$ is exact.

$\Box$ \vskip 1cm

We recall
the \emph{proximal relation $P \subset X \times X$} in $(X,f)$. We say that $(x, y) \in P$ if $\liminf d(f^n x, f^n y) = 0,  \ n \in \mathbb{N}$. Equivalently, for any $\epsilon > 0$ there is a $n \in \mathbb{N}$ such that $d(f^n x, f^n y) < \epsilon$. The relation $P$ is reflexive, symmetric, and $f$ invariant.The system $(X, f)$ is said to be \emph{distal} if \textbf{$P = \bigtriangleup$}, i.e. there are no non-trivial proximal pairs. It is obvious that in such a case our map will be a homeomorphism.

The system  $(X, f)$ is \emph{equicontinuous} if the maps of $X$ defined by the iterates of $f$ form an equicontinuous family. That
is, for every $\epsilon > 0$, there is a $\delta > 0$ such that whenever $d(x, x') < \delta$ then $d(f^n x, f^n x') < \epsilon$
for all $n \in \mathbb{N}$.

We know that equicontinuous systems are distal, and it is well known that the converse fails. We note that the system $(X,f)$ is equicontinuous if and only if the system $(2^X, f_*)$ is equicontinuous \cite{sn1}. However, we observe that $(X,f)$ distal need not imply that $(2^X, f_*)$ is distal.

\vskip .2cm

For the unit disc $D$, consider the map $f : D \rightarrow D$  defined by $f(r, \theta) = (r, r + \theta)$. Clearly, $(D,f)$ is distal.
Let $r_n \mapsto i$, where $r_n$ are rationals and $i$ an irrational, and consider $A = \{ (r_n, 0) : n \in \mathbb{N} \} \cup \{ (i,0) \} \in 2^X$.   The orbit closure of $A$ will consist of the set $B$ with the full
circle at the irrational radius $i$ and  periodic orbits at each rational
radius $r_n$. So $A$ is not in the orbit closure of $B$.

A point is called \emph{almost periodic} if its orbit closure is a minimal system. It is known that a distal system can be realized as a union of minimal systems, and so every point is  an \emph{almost periodic point}. Since a distal system is pointwise almost periodic, and the point $A$ fails to be almost periodic,  $(2^D, f_*)$ cannot be distal.

To this end, we recall the definition of $RP \subset X \times X$ the \emph{regionally proximal relation}. We say that $(x,y) \in RP$ if there are sequence $\{x_n\}$ and $\{y_n\}$ with $x_n \to x$ and $y_n \to y$ and integers $k_n$ such that $d(f^{k_n}(x_n),f^{k_n}(y_n))\to 0$. It is  known  that $(X,f)$ is equicontinuous if and only if $RP=\Delta$. That is, whenever $x \neq y$, $(x,y)\notin RP$.

 A classical theorem of Ellis says that $(X,f)$ is distal if and only if the product flow $(X \times X, f \times f)$ is pointwise almost periodic. That is, every orbit closure in $X \times X$ is minimal.

We now prove when $(2^X, f_*)$ will be distal.

\begin{theo} The following are equivalent:

$(i)$ $(X,f)$ is equicontinuous.

$(ii)$ $(2^X,f)$ is equicontinuous

$(iiii)$ $(2^X,f)$ is distal.

\end{theo}

\proof It is easy to show that $(i) \implies (ii)$, and as since equicontinuous systems are distal, $(ii) \implies (iii)$.

For $(iii) \implies (i)$, suppose that $(X,f)$ is not equicontinuous.
Then there are $x$ and $y$ in $X$ with $x \neq y$ and $(x,y) \in RP$. Then we have $x_n \to x$, $y_n \to y$ and $k_n \in \Bbb Z$ with $(f^{k_n}(x_n),f^{k_n}(y_n))\to (z,z)$ for some $z \in X$. Let $C=\{x_1,x_2, \dots,x\}$ and $D=\{y_1,y_2, \dots, y\}$. We may suppose $C \cap D=\emptyset$. We show that the orbit closure of $(C,D)$ in $2^X \times 2^X$ is not minimal. Let (a subsequence of) $(f_* \times f_*)^{k_n}(C ,D) \to (C',D')$. Then $C' \cap D' \neq \emptyset$ (in the terminology of Proposition 1.3, $(C',D') \in INT$). Now if the orbit closure of $(C,D)$ were minimal, there would be a sequence $\{l_m\}$ with  $(f_* \times f_*)^{l_m}(C',D') \to (C,D)$. But in this case $C \cap D \neq \emptyset$, a contradiction.

$\Box$

A system $(X,f)$ is  \emph{sensitive} if there exists a $\delta > 0$ such that for any $x_1 \in X$ and opene set $U \ni x_1$, there exists $ x_2  ( \neq x_1) \in U$ and $m \in \mathbb{N}$ such that $d(f^m(x_1), f^m(x_2)) > \delta$.

It is known that minimal systems are either equicontinuous or sensitive.  Various forms of sensitivities for the induced system $(2^X,f_*)$ are studied by Sharma and Nagar in \cite{sn}, and among other things there is an example of a sensitive $(X,f)$ for which $(2^X,f_*)$ is not sensitive.

%eja I changed the title since the previous section also concerned
\section{ $(X,f)$ Can Be Isomorphic to Its Own Induced System }

In this section we will again restrict attention to compact systems $(X,f)$.

 The dynamical system $(2^X,f_*)$ induces the system $(2^{2^X}, f_{**})$, and in many cases we believe that this will be a transfinite process. Hence, it would be nice to characterize those compact dynamical systems $(Y,g)$ which are of the form
$(2^X,f_*)$ for some $(X,f)$, (necessarily a compact system since $X$ is isometric with a closed subset of $2^X$). In this section, we look into some cases where this process of  inducing  definitely gives us nonconjugate dynamics  or where this process of inducing stabilizes.

 A subset $A \subset X$
is called \emph{+ invariant} when $f(A) \subset A$ and invariant when $f(A) = A$. Observe that if $A$ is closed
then it is invariant if and only if as an element of $2^X$ it is a fixed point for $f_*$. Hence, if $f$ is transitive then
if $A$ is a fixed point for $f_*$ then either $A = X$ or $A$ has empty interior. If $f$ is \emph{totally transitive}, i.e. every
power $f^n$ for $n = 1,...$ is transitive, then if $A$ is a periodic point for $f_*$, and hence a fixed point for some
$f^n_*$, then either $A = X$ or $A$ has empty interior. In particular, a weak mixing system is totally transitive and so this then holds.

If $F \subset 2^X$ is a closed invariant set for $f_*$ then
as an element of $2^{2^X}$ it is a fixed point and hence $\vee(F) = \bigcup F$ is a fixed point for $f_*$, i.e. an
invariant set for $f$.  In particular the union of any periodic orbit of $f_*$ is a fixed point of $f_*$. Notice that the
union of any finite collection of periodic points of $f_*$ is a periodic point (with period dividing the l.c.m. of the
periods of the members of the collection). The intersection, if non-empty, is also a periodic point.

Recall that if the periodic points are dense in $X$ then the same holds for $2^X$.  The converse is not true as we will see below.
%As a first example, let $(X_1,f_1)$ be the shift map on $X_1 = \{0,1\}^{\Z}$ and $(X_2,f_2)$ be a weak mixing minimal system.
%The product system $(X_1 \times X_2, f_1 \times f_2)$ is weak mixing. Since it has a nontrivial minimal factor it has no periodic points, but each $\{p\}\times X_2$ is periodic when $p$ is.
Any odometer $(X,f)$ is minimal, and equicontinuous and (so is not totally transitive). Any point of $x$ can be approximated
by a periodic element of $2^X$ with nonempty interior and so with the union of the orbit equal to $X$. It follows that the
periodic points are dense in $2^X$.

\begin{theo}\label{theo4.1} Let $(X,f)$ be a compact dynamical system. The periodic points are dense in $(2^X,f_*)$ if their
closure contains $i(X) \subset 2^X$.

If $f$ is totally transitive and the periodic points are dense in $2^X$ then there are infinitely many distinct fixed points
for $f_*$.\end{theo}

\proof If $\{x_1,..,x_k\} \subset X$ and $\ep >0 $ then there exist periodic points $\{A_1,..., A_k \} \subset 2^X$
such that $d(A_i, \{x_i\}) < \ep $ for $i = 1,...,k$. Then $A = A_1 \cup \ldots \cup A_k$ is a periodic point with
$d(A, \{x_1, \ldots ,x_k\} ) < \ep $.  Hence, the closure of the periodic points contains $FIN(X)$ which is dense in $2^X$.

The unions of each of a finite collection of periodic orbits for $f_*$, other than the fixed point $X$ itself, are
closed, nowhere dense sets since $f$ is totally transitive. Each of the unions is a fixed point. There is an opene set
$U \subset X$ disjoint from all of them and it contains a periodic point $A$ for $f_*$.  The union of the associated periodic
orbit is not contained in the previous union and so is distinct from the previous fixed points.

$\Box$ \vspace{.5cm}

\begin{cor}\label{cor4.2} If $(X,f)$ is a mixing subshift of finite type then it is not isomorphic to $(2^Y,g_*)$
for any compact system $(Y,g)$.\end{cor}

\proof A subshift of finite type has dense periodic points and is expansive. Hence, it cannot admit infinitely many
distinct fixed points.

$\Box$ \vspace{.5cm}

An \emph{inverse system} is a sequence between compact metrizable spaces: $\{ \pi_i : X_{i+1} \to X_{i} : i = 1,2,...\}$.
The associated \emph{inverse limit} is the compact metrizable space $X_{\infty} = \{ x \in \Pi_i \ X_i : x_i = \pi_i(x_{i+1}) $ for $ i = 1,2,... \}$, which is metrizable when the $X_i$'s are.
The projection map the restriction map to the $i^{th}$ coordinate is denoted $\pi^i : X_{\infty} \to X_{i}$. If $A$ is a closed
subset of $X_{\infty}$ we let $F_i = \pi^{i}(F)$. If $x \in X_{\infty}$ is such that $\pi_i(x) \in F_{i}$ for all $i$ then
$\{ F \cap (\pi_i)^{-1}(x_i) \}$ is a decreasing sequence of nonempty compacta and so the intersection, which is
$F \cap \{ x \}$ is nonempty. That is, $x \in F$. Thus, the restrictions $\pi^i : F \to F_i$ and $\pi_i : F_{i+1} \to F_i$ are
all surjective and so $F$ is the inverse limit of the system $\{ \pi_i : F_{i+1} \to F_i \}$. It thus follows that
the map from $\{ (\pi_i)_* : 2^{X_{i+1}} \to 2^{X_i} \} $ is an inverse system.  Furthermore the
maps $(\pi^i)_* : 2^{X_{\infty}} \to 2^{X_i}$ induce a natural homeomorphism between $2^{X_{\infty}}$
and the inverse limit of the system
$\{ (\pi_i)_* : 2^{X_{i+1}} \to 2^{X_i} \} $. Refer to \cite{in, se} for more details.

Now let $X$ be any compact metric space. Define $X_1 = 2^X, X_2 = 2^{2^X} = 2^{X_1}$ and $\pi_1 = \vee : X_2 \to X_1$.
Inductively, define $X_{i+1} = 2^{X_i}$ and $\pi_{i} = (\pi_{i-1})_* : X_{i+1} = 2^{X_i} \to 2^{X_{i-1}} = X_i$. Let
$X_{\infty}$ denote the inverse limit of this inverse system. As remarked in the previous paragraph, there is a
natural homeomorphism from $2^{X_{\infty}}$ to the inverse limit of the system $\{ (\pi_i)_* : 2^{X_{i+1}} \to 2^{X_i} \}$
which is our
original system with the index shifted by one. That is, there is a natural homeomorphism  $q :2^{X_{\infty}} \to X_{\infty}$.

If $(X,f)$ is a compact dynamical system then letting $f_1 = f_*$ and $f_{i+1} = (f_i)_*$ we obtain an inverse limit of dynamical
systems with $f_{\infty}$ the map induced on $X_{\infty}$ and $q$ is an isomorphism from $(2^{X_{\infty}},(f_{\infty})_*)$ to
$(X_{\infty},f_{\infty})$. Since the inverse limit of weak mixing systems is weak mixing, $(X_{\infty},f_{\infty})$ is
weak mixing when $(X,f)$ is. Furthermore, the pair of action maps
$(\pi_1)_* : 2^{X_{\infty}} \to 2^{2^X}$ and $\pi_1 : X_{\infty} \to 2^X$ maps the isomorphism $q$ to $\vee$.  Thus, we have proved
the following result.

\begin{theo}\label{4.3} If $(X,f)$ is a compact dynamical system then there $(2^X,f_*)$ is a factor of
a compact dynamical system $(X_{\infty},f_{\infty})$ which is isomorphic to $(2^{X_{\infty}},(f_{\infty})_*)$. Furthermore, if
$(X,f)$ is weak mixing then $(X_{\infty},f_{\infty})$ is weak mixing. \end{theo}

$\Box$ \vspace{.5cm}

If $\{ n_i \}$ is a strictly increasing sequence of positive integers with $n_i | n_{i+1}$ then
the inclusions of $n_{i+1}\Z \to n_i\Z$ induces an inverse system $\{ \pi_i : \Z/n_{i+1}\Z \to \Z/n_i\Z \}$ of epimorphisms
of finite rings.  The inverse limit $\Z_{\{n_i\}}$ is a compact ring which is topologically a Cantor set. The associated
\emph{odometer} or \emph{adding machine} uses the translation map $t(x) = x + 1$. Every clopen subset is a periodic
element for $(2^{\Z_{\{n_i\}}},t_*)$, as it is the preimage of a finite subset of $\Z/\Z_{\{n_i\}}$ for sufficiently
large $i$. Hence, from Theorem \ref{theo4.1} it follows that the periodic points are dense in
$2^{\Z_{\{n_i\}}}$.  On the other hand, $(\Z_{\{n_i\}},t)$ is minimal and so has no periodic points.

An interesting question is which dynamical systems can be represented as $(2^X,f_*)$ for some system $(X,f)$. This is far from settled, but we have obtained some obstructions to such. Of course, $2^X$ can never be minimal. Theorem  5.3 tells us that it can't be backwards minimal unless it is exact, and by Theorem 5.4, it can't be distal unless it is equicontinuous.

Moreover, $(2^X,f_*)$ has some structure which apparently does not occur in most dynamical systems, namely the relations $INT, INC$, and $EPS$ defined in Proposition 1.3. These are closed $f_* \times f_*$ invariant relations, and an interesting problem is to determine their dynamical significance.

\vskip .5cm

\textbf{Acknowledgement}: This work was done when the third author visited University of Maryland. She acknowledges the hospitality of the Mathematics Department at the University.

\bibliography{xbib}

\end{document}